\documentclass[12pt]{amsart}
\usepackage{amsmath}
\usepackage{amsfonts}
\usepackage{amssymb}
\usepackage{amsthm}
\usepackage{amstext,amscd}

\theoremstyle{plain}
\newtheorem{theorem}{Theorem}[section]

\newtheorem{conjecture}{Conjecture}
\newtheorem{corollary}[theorem]{Corollary}

\newtheorem{definition}[theorem]{Definition}

\newtheorem{lemma}[theorem]{Lemma}

\newtheorem{proposition}[theorem]{Proposition}
\newtheorem{remark}{Remark}

\numberwithin{equation}{section}
\begin{document}
\title{An $SL_2(\mathbb{C})$ Algebro-Geometric Invariant of Knots}
\author{Weiping Li}
\address{Department of Mathematics\\
         Oklahoma State University\\
         Stillwater, OK 74078}
\email{wli@math.okstate.edu}
\author{Qingxue Wang}
\address{Chern Institute of Mathematics\\
         Nankai University\\
         Tianjin 300071\\
         P.R. China}
\email{qingxue\_wang@yahoo.com.cn}

\begin{abstract}
In this paper, we define a new algebro-geometric invariant of
$3$-manifolds resulting from the Dehn surgery along a hyperbolic
knot complement in $S^3$. We establish a Casson type invariant for
these $3$-manifolds. In the last section, we explicitly calculate
the character variety of the figure-eight knot and discuss some
applications, as well as the computation of our new invariants for
some 3-manifolds resulting from the Dehn surgery along the
figure-eight knot.

\end{abstract}
\maketitle

\section{Introduction}
It is well-known that there are successful applications for the
$SU(2)$-Casson invariants of knots in integral homology 3-spheres,
so it is a natural question to extend them to other compact or
noncompact groups. For a hyperbolic 3-manifold, by Thurston's result
(\cite{Th}), there is a standard (up to conjugacy) $SL(2,
\mathbb{C})$-representation associated to the hyperbolic structure.
One would like to have the $SL(2, \mathbb{C})$-Casson invariant for
this representation. There are some attempts in this direction, see
\cite{Cu1, Cu2, Li}. But the difficulties lie in the fact that not
so much information is known about the character variety of a knot
group. In this paper, we shall approach this problem purely from the
algebraic geometric point of view.

Consider a hyperbolic knot $K$ in $S^3$. Let $X_0$ be the
irreducible component of the character variety of $K$ which contains
the character of the discrete faithful representation associated to
the hyperbolic structure of $K$. It is known that $X_0$ is an affine
curve. We study in detail the properties of its image in the
character variety of the boundary (which is a torus) and the
corresponding part in the $A$-polynomial. Using $X_0$, we construct
a new $SL_2(\mathbb{C})$ algebro-geometric invariant $\lambda(p,q)$
(see Definition \ref{lambda}) for the manifolds obtained by the
$(p,q)$ Dehn surgery along the knot complement. Roughly speaking,
our invariant $\lambda(p,q)$ is counting the geometric intersection
multiplicity of $X_0$ with another affine curve in the character
variety of the boundary. It does not require the non-sufficiently
large condition in \cite{Cu1}. The non-vanishing of $\lambda (p, q)$
implies that the 3--manifold $K(p/q)$ resulting from $(p, q)$
surgery is non-cyclic (Proposition \ref{prop-liv}). Via our
invariant and Culler-Shalen norm, we obtain an upper bound for the
number of ideal points with multiplicity which are zeroes of the
function $f_{\gamma}$ (Proposition \ref{deg-ine}). The definition of
$f_{\gamma}$ is in Section $3$ and its degree is the Culler-Shalen
norm.

Our point of view is that the component $X_0$ should contain a lot
of topological and geometric information about the knot and its Dehn
fillings. Of course, the whole character variety may contain more
information, but at present we still have little information about
other components. It seems that the whole character variety is
complicated in general. That is why we restrict ourselves to the
single component $X_0$ in this paper.

The paper is organized as follows. In section $2$, we introduce the
notations used in the paper. In section $3$, we define the
algebro-geometric invariant $\lambda(p,q)$ of a hyperbolic knot, and
study its properties. In the last section, we explicitly calculate
several invariants of the figure-eight knot, for instance, its
character variety, ideal point, Culler-Shalen norm and our new
invariant. Some applications of these calculations are also
discussed.

\section{Terminology and Notation}
\subsection{}
Let $K$ be a knot in $S^3$ and $M_K$ its complement. That is,
$M_K=S^3-N_K$ where $N_K$ is the open tubular neighborhood of $K$ in
$S^3$. $M_K$ is a compact $3$-manifold with boundary $\partial
M_K=T^2$ a torus. Denote by
$R(M_K)=\text{Hom}(\pi_1(M_K),SL_2(\mathbb{C}))$ and $R(\partial
M_K)=\text{Hom}(\pi_1(\partial M_K),SL_2(\mathbb{C}))$. It is known
that they are affine algebraic sets over the complex numbers
$\mathbb{C}$ and so are the corresponding character varieties
$X(M_K)$ and $X(\partial M_K)$ (See \cite{CS}). We also have the
canonical surjective morphisms $t:R(M_K)\longrightarrow X(M_K)$ and
$t:R(\partial M_K)\longrightarrow X(\partial M_K)$ which map a
representation to its character. The natural homomorphism $i:
\pi_1(\partial M_K)\longrightarrow \pi_1(M_K)$ induces the
restriction maps $r: X(M_K)\longrightarrow X(\partial M_K)$ and $r:
R(M_K)\longrightarrow R(\partial M_K)$.
\subsection{}
Throughout this paper, for a matrix $A\in SL_2(\mathbb{C})$, we
denote by $\sigma(A)$ its trace.

\subsection{}
Since $\pi_1(\partial M_K)=\mathbb{Z}\oplus \mathbb{Z}$, we shall
fix two oriented simple curves $\mu$ and $\lambda$ as its
generators. They are called the meridian and longitude respectively.
Let $R_D$ be the subvariety of $R(\partial M_K)$ consisting of the
diagonal representations. Then $R_D$ is isomorphic to
$\mathbb{C}^{*} \times \mathbb{C}^{*}$. Indeed, for $\rho \in R_D$,
we obtain
\begin{equation*}
\rho(\mu)=
  \left[
  \begin{matrix}
    m & 0\\
    0 & m^{-1}
  \end{matrix}
  \right], \;\;
\rho(\lambda)=
  \left[
  \begin{matrix}
    l & 0\\
    0 & l^{-1}
  \end{matrix}
  \right].
\end{equation*}
then we assign the pair $(m,l)$ to $\rho$. Clearly this is an
isomorphism. We shall denote by $t_{D}$ the restriction of the
morphism $t:R(\partial M_K)\longrightarrow X(\partial M_K)$ on
$R_D$.

\section{$SL_2(\mathbb{C})$ Character Variety and a New Knot
Invariant} In this section, let $K$ be a hyperbolic knot in $S^3$.
Then $M_K$ is a hyperbolic $3$-manifold with finite volume. Since
$M_K$ is hyperbolic, there is a discrete faithful representation
$\rho_0\in R(M_K)$ corresponding to its hyperbolic structure. We
denote by $R_0$ an irreducible component of $R(M_K)$ containing
$\rho_0$. Let $X_0=t(R_0)$. By \cite[Proposition~1.1.1]{CGLS},
$X_0\subset X(M_K)$ is an irreducible affine variety of dimension
$1$.

\subsection{The curve of characters} In this subsection, we give some
elementary properties of the character varieties and its component
$D_0$.

Since $\partial M_K$ is a torus, we identify $R(\partial M_K)$ with
the set $\{(A,B)|A,B\in SL_2(\mathbb{C}), AB=BA\}$, where
$A=\rho(\mu)$, $B=\rho(\lambda)$ for $\rho$ a representation of
$\pi_1(\partial M_K)$ in $SL_2(\mathbb{C})$. As in Section $2.3$,
$R_{D}\subset R(\partial M_K)$ is the subvariety consisting of the
representations of diagonal matrices. We have the isomorphism
$p:R_{D}\rightarrow \mathbb{C}^{*}\times \mathbb{C}^{*}$ defined by
$p(\rho)=(m,l)$ if
\begin{equation*}
\rho(\mu)=
  \left[
  \begin{matrix}
    m & 0\\
    0 & m^{-1}
  \end{matrix}
  \right], \;\;
\rho(\lambda)=
  \left[
  \begin{matrix}
    l & 0\\
    0 & l^{-1}
  \end{matrix}
  \right].
\end{equation*}

Let $R_{D}$ be identified with $\mathbb{C}^{*}\times \mathbb{C}^{*}$
via $p$. By the proof of \cite[Proposition~1.4.1]{CS}, $\chi\in
X(\partial M_K)$ is determined by its values on $\mu$, $\lambda$ and
$\mu\lambda$. Define a map $t:R(\partial M)\rightarrow \mathbb{C}^3$
by
$t(\rho)=(\sigma(\rho(\mu)),\sigma(\rho(\lambda)),\sigma(\rho(\mu\lambda)))$
. Then $X(\partial M_K)=t(R(\partial M_K))$. This map is the regular
surjective morphism $t:R(\partial M_K)\rightarrow X(\partial M_K)$.
That is why we use the same letter $t$. The map $t_D$, the
restriction of $t$ on $R_D=\mathbb{C}^{*}\times \mathbb{C}^{*}$, is
given explicitly by, for $(m,l)\in \mathbb{C}^{*}\times
\mathbb{C}^{*}$,
\[
t_D(m,l)=(m+m^{-1}, l+l^{-1}, ml+m^{-1}l^{-1}).
\]
It is straightforward to check that $t_D(\mathbb{C}^{*}\times
\mathbb{C}^{*})=X(\partial M_K)$ and $t_D$ is $2:1$ except at four
points $(\pm 1, \pm 1)$ where it is $1:1$.

We have the following diagram:

\[
  \begin{CD}
    R(M_K)\supset R_0   @. D_0=t^{-1}_{D}(Y_0)\subset \mathbb{C}^{*}\times \mathbb{C}^{*}=R_D\subset R(\partial M_K)\\
     @VtVV   @Vt_DVV\\
    X(M_K)\supset X_0=t(R_0) @>r>> Y_0=\overline{r(X_0)}\subset X(\partial
    M_K).
\end{CD}
\]

First we characterize the map $t_D$ and the character variety
$X(\partial M)$.
\begin{proposition}\label{f-prop}
The map $t_D:R_D=\mathbb{C}^{*}\times \mathbb{C}^{*}\rightarrow
X(\partial M_K)$ is a finite morphism.
\end{proposition}

\begin{proof}
The affine coordinate ring for $\mathbb{C}^{*}\times \mathbb{C}^{*}$
is $\mathbb{C}[x,x^{-1},y,y^{-1}]$. Notice that $t_D$ is surjective,
we can identify the coordinate ring of $X(\partial M_K)$ with the
sub-ring $\mathbb{C}[x+x^{-1},y+y^{-1},xy+x^{-1}y^{-1}]$ of
$\mathbb{C}[x,x^{-1},y,y^{-1}]$. Let $t=x+x^{-1}$, then $x$ and
$x^{-1}$ are roots of the equation $X^2-tX+1=0$. Hence $x$ and
$x^{-1}$ are integral over
$\mathbb{C}[x+x^{-1},y+y^{-1},xy+x^{-1}y^{-1}]$, so are $y$ and
$y^{-1}$. Now $x$, $x^{-1}$, $y$, $y^{-1}$ are integral over
$\mathbb{C}[x+x^{-1},y+y^{-1},xy+x^{-1}y^{-1}]$, it follows that
$\mathbb{C}[x,x^{-1},y,y^{-1}]$ is also integral over it. Therefore,
$t_D$ is finite.
\end{proof}

\begin{proposition}\label{torus}
The character variety $X(\partial M_K)$ is a surface in
$\mathbb{C}^3$ defined by
  \begin{equation}\label{cha-torus}
    x^2+y^2+z^2-xyz-4=0.
  \end{equation}
\end{proposition}

\begin{proof}
Let $x$, $y$, $z$ be the traces of $a=\rho(\mu)$, $b=\rho(\lambda)$,
$c=\rho(\mu\lambda)$ respectively. By the formula (\ref{tf4}), we
have
\[
  \sigma((ab)a^{-1}b^{-1})=\sigma(ab)\sigma(a^{-1}b^{-1})+\sigma(a^{-1})\sigma(a)+\sigma(b^{-1})\sigma(aba^{-1})
  -\sigma(ab)\sigma(a^{-1})\sigma(b^{-1})-\sigma(I).
\]

Since $a$ and $b$ commute, $\sigma((ab)a^{-1}b^{-1})=\sigma(I)=2$,
where $I$ is the $2\times 2$ identity matrix. This gives the
equation (\ref{cha-torus}) by (\ref{tf1}) and (\ref{tf2}). Hence,
$X(\partial M_K)$ is contained in the surface in $\mathbb{C}^3$.

On the other hand, let $(x,y,z)\in \mathbb{C}^3$ be a solution to
(\ref{cha-torus}). It is a straightforward calculation that we can
find $(m,l)\in \mathbb{C}^{*}\times \mathbb{C}^{*}$ such that
$x=m+m^{-1}$, $y=l+l^{-1}$, and $z= ml+m^{-1}l^{-1}$. Hence the
result follows.
\end{proof}

For each $\gamma \in \pi_1(M_K)$, there is a natural regular map
$I_{\gamma}: X(M_K)\rightarrow \mathbb{C}$ defined by
$I_{\gamma}(\chi)=\chi (\gamma)$. The set of functions $I_{\gamma}$,
$\gamma \in \pi_1(M_K)$ generates the affine coordinate ring of
$X(M_K)$. Moreover, by \cite[Proposition~1.1.1]{CGLS}, for each
nonzero $\gamma \in \pi_1(\partial M_K)$, the function $I_{\gamma}$
is non-constant on $X_0$. This implies that $r(X_0)\subset
X(\partial M_K)$ has dimension $1$. We set $Y_0=\overline{r(X_0)}$,
the Zariski closure of $r(X_0)$ in $X(\partial M_K)$. Then $Y_0$ is
an irreducible affine curve. Denote by $D_0$ the inverse image
$t^{-1}_{D}(Y_0)$.

\begin{proposition}\label{d-curve}
(i) The inverse image $D_0\subset \mathbb{C}^{*}\times
\mathbb{C}^{*}$ is an affine algebraic set of dimension $1$;\\
(ii) the image of each $1$-dimensional component of $D_0$ under
$t_D$ is the whole $Y_0$;\\
(iii) $D_0$ has no $0$-dimensional components and has at most two
$1$-dimensional components.
\end{proposition}

\begin{proof}
(i). As $Y_0$ is closed, so is $D_0=t^{-1}_{D}(Y_0)$. Thus $D_0$ and
$Y_0$ have the same dimension $1$ because $t_D$ is finite by
Proposition \ref{f-prop}.

(ii). Next we show that if $V$ is a $1$-dimensional irreducible
component of $D_0$, then $t_D(V)=Y_0$. In fact, since a finite
morphism maps closed sets to closed sets, $t_D(V)\subset Y_0$ is
closed, irreducible and its dimension is $1$. Hence $t_D(V)=Y_0$ due
to the irreducibility of $Y_0$.

(iii). Suppose that $D_0$ has three distinct $1$-dimensional
components $V_i$, $1\leq i \leq 3$. Since the $V_i \cap V_j$, $i\ne
j$ are empty or finite sets, we can choose $y_0\in Y_0$ such that
$y_0\notin t_D(V_i \cap V_j)$ for $\forall i\ne j$. We see that
$t_D(V_i)=Y_0$, hence $t_D^{-1}(y_0)$ has three elements. This
contradicts that $t_D^{-1}(y)$ has at most two elements for any
$y\in Y_0$. Therefore, $D_0$ has at most two $1$-dimensional
irreducible components.

Now we show that it has no $0$-dimensional components. Consider the
morphism $\tau: \mathbb{C}^{*}\times \mathbb{C}^{*}\rightarrow
\mathbb{C}^{*}\times \mathbb{C}^{*}$ defined by
$\tau(m,l)=(m^{-1},l^{-1})$. It is an involution and hence an
automorphism of $\mathbb{C}^{*}\times \mathbb{C}^{*}$ with
$\tau(D_0)=D_0$. So it is an automorphism of $D_0$ when restricting
on it. Therefore it maps $1$-dimensional component to
$1$-dimensional component. Now for any $y\in Y_0$, $\tau$ permutes
elements of $t_D^{-1}(y)$ and at least one element of $t_D^{-1}(y)$
is contained in some $1$-dimensional component since $t_D$ maps
every $1$-dimensional component onto $Y_0$. Therefore, no
$0$-dimensional component for $D_0$.
\end{proof}

\begin{remark}
For the character variety $X(M_K)$, by
\cite[Proposition~2.4]{CCGLS}, there is no zero-dimensional
components of $X(M_K)$. Proposition \ref{d-curve} part (iii) shows
that $D_0$ in $R_D\subset R(\partial M_K)$ has no $0$-dimensional
components.
\end{remark}

In summary, we have the following two cases:\\
(I) $D_0$ itself is an irreducible affine curve;\\
(II) $D_0=V_1\cup V_2$, where $V_i$ are two $1$-dimensional
irreducible components. Each is an irreducible affine curve with
$t_D(V_i)=Y_0$, $i=1,2$. Moreover, they are isomorphic to each other
under the involution $\tau$ by Proposition \ref{d-curve} (iii) . We
also have $t_D: V_i\rightarrow Y_0$ is a one-to-one and onto regular
map.

\subsection{Invariant and Dehn Surgery of Hyperbolic Knots}
In this subsection, we define an algebro-geometric invariant of the
$3$-manifolds $K(p/q)$ resulting from the $(p,q)$ Dehn surgery along
a hyperbolic knot complement in $S^3$. A Casson type
$SL_2(\mathbb{C})$ invariant for $K(p/q)$ is also established.

Let $A_0(m,l)$ be the defining equation of the closure of the affine
curve $D_0$ in $\mathbb{C}\times \mathbb{C}$. We also require that
it has no repeated factors. It is a factor of the $A$-polynomial of
the knot $K$ defined in \cite{CCGLS}, and $A_0(m,l)^d$ is the
$A$-polynomial of $X_0$ defined in \cite[Page~109]{BZ}, where $d$ is
the degree of regular map $r:X_0\rightarrow Y_0$.

Denote by $\widetilde{X_0}$ (resp. $\widetilde{Y_0}$) a smooth
projective model of the affine curve $X_0$ (resp. $Y_0$). The
restriction morphism $r:X_0\rightarrow Y_0$ induces a regular map
$\widetilde{r}: \widetilde{X_0}\rightarrow \widetilde{Y_0}$.

\begin{lemma}\label{le1}
The regular map $\widetilde{r}$ is an isomorphism.
\end{lemma}

\begin{proof}
By our assumption, $M_K$ is the complement of a knot in $S^3$, hence
$H^1(M_K;\mathbb{Z}_2)=\mathbb{Z}_2$. By \cite[Corollary~3.2]{Dun},
$r$ is a birational isomorphism onto $r(X_0)$. Since $\widetilde{r}$
is induced by $r$, the result follows.
\end{proof}

Let $\gamma=p\mu+q\lambda\in H_1(\partial M_K;\mathbb{Z})$ be a
non-zero primitive element with $p,q$ coprime. Define a regular
function $f_{\gamma}=I_{\gamma}^2-4$ on $X_0$. Since $I_{\gamma}$ is
nonconstant, so is $f_{\gamma}$. It is also a meromorphic function
on $\widetilde{X_0}$ or equivalently, a non-constant holomorphic
function from $\widetilde{X_0}$ to $\mathbb{CP}^1$ and we denote it
again by $f_{\gamma}$ .

Since $\gamma\in H_1(\partial M_K;\mathbb{Z})$, we can think of
$I_{\gamma}$ as a regular function on $Y_0$ in $X(\partial M)$.
Define on $Y_0$ the function $f_{\gamma}^{\prime}=I_{\gamma}^2-4$.
Similarly, it is a non-constant regular function on $Y_0$, and hence
a non-constant holomorphic function from $\widetilde{Y_0}$ to
$\mathbb{CP}^1$, denoted also by $f_{\gamma}^{\prime}$. Then by the
definition, we have
\[
f_{\gamma}=f_{\gamma}^{\prime}\circ \widetilde{r}.
\]
In particular, $\text{deg}f_{\gamma}=\text{deg}f_{\gamma}^{\prime}$
by Lemma \ref{le1}.

We denote by $Z_{\gamma}$ the set of zeroes of the function
$f_{\gamma}$ on $X_0$. If $\chi\in Z_{\gamma}$, then there exists a
representation $\rho\in R_0$ such that its character
$\chi_\rho=\chi$.

\begin{lemma}\label{le2}
Suppose $\chi_\rho\in Z_{\gamma}$. Then either $\rho(\gamma)=\pm
\text{I}$ or $\rho(\gamma)\ne \pm I$ and $\sigma(\rho(\alpha))=\pm
2$ for all $\alpha\in \pi_1(\partial M_K)$.
\end{lemma}

\begin{proof}
Note that $f_{\gamma}(\chi_\rho)=0$ is equivalent to that the trace
of $\rho(\gamma)$ is $\pm 2$. If $\rho(\gamma)\ne \pm I$, then
$\rho(\gamma)$ is a parabolic element in $SL_2(\mathbb{C})$. Since
$\rho(\gamma)$ and $\rho(\alpha)$ commute for all $\alpha\in
\pi_1(\partial M_K)$, $\sigma(\rho(\alpha))=\pm 2$.
\end{proof}

Suppose that $\rho\in R_0\subset R(M_K)$ is irreducible. Then its
character $\chi_{\rho}$ is contained in the component $X_0$. Assume
that $\rho(\mu)$ and $\rho(\lambda)$ are parabolic. Up to
conjugation, we have
\begin{equation*}
\rho(\mu)=\pm
  \left[
  \begin{matrix}
    1 & 1\\
    0 & 1
  \end{matrix}
  \right]; \;\;
\rho(\lambda)=\pm
  \left[
  \begin{matrix}
    1 & t(\rho)\\
    0 & 1
  \end{matrix}
  \right],
\end{equation*}

where $t(\rho)$ is a complex number. We have the following:
\begin{conjecture}
Let $\rho\in R_0$ be an irreducible
$SL_2(\mathbb{C})$-representation of a hyperbolic knot in $S^3$. If
$\rho(\mu)$ and $\rho(\lambda)$ are parabolic, then $t(\rho)\notin
\mathbb{Q}$.
\end{conjecture}

\begin{remark}
(1) If $\rho=\rho_0$ is the discrete faithful representation of the
hyperbolic structure, then we know $\rho(\mu)$ and $\rho(\lambda)$
are parabolic and $t(\rho_0)$ is called the \emph{cusp constant} and
the \emph{cusp polynomial} is the minimum polynomial for $t(\rho_0)$
over $\mathbb{Q}$. Moreover $\{1, t(\rho_0)\}$ generates a lattice
of $\mathbb{C}$. Therefore, $t(\rho_0)\notin \mathbb{R}$. For more
detail, see \cite[Section~6]{CL}. So $t(\rho)$ can be thought of as
the generalization of the cusp constant.

(2) The conjecture can not be extended to non-hyperbolic knots.
There is an example in \cite[Section~5]{Ril1} of an irreducible
parabolic representation $\rho$ of alternating torus knot such that
$t(\rho)\in \mathbb{Z}$. It is unlikely that the conjecture can be
extended to the non-geometric components of a hyperbolic knot in
$S^3$, that is, components not containing any representation of the
hyperbolic structure of the knot.

(3) The conjecture fails for knots in other rational homology
spheres. For instance, for the SnapPea manifold m023, which is the
complement of a knot in the lens space $L(3,1)$, the geometric
component contains a representation with $t(\rho)\in \mathbb{Q}$.
\end{remark}

The conjecture is true for the figure-eight knot. In fact, for the
figure-eight knot, by Proposition \ref{8cha-1}, if $\rho$ satisfies
the condition of the conjecture, then it must be the discrete
faithful representation of the hyperbolic structure. By the remark,
$t(\rho)\notin \mathbb{R}$.

\begin{proposition}\label{kin}
Suppose the above conjecture is true. Let $\gamma=p\mu+q\lambda\in
H_1(\partial M_K;\mathbb{Z})$ with $p,q$ non-zero coprime integers.
Let $\chi_{\rho}\in Z_{\gamma}$ be the character of an irreducible
representation $\rho$. Then $\rho(\gamma)=\pm I$ if and only if the
trace of $\rho(\mu)$ is not equal to $\pm 2$.
\end{proposition}

\begin{proof}
If the trace of $\rho(\mu)$ is not equal to $\pm 2$, then
$\rho(\mu)$ is not parabolic or $\pm I$. Since $\rho(\mu)$ and
$\rho(\gamma)$ commute, $\rho(\gamma)$ is not parabolic. Since
$\chi_{\rho}\in Z_{\gamma}$, $\sigma(\rho(\gamma))=\pm 2$. Thus,
we obtain $\rho(\gamma)=\pm I$.\\

Next suppose that the trace of $\rho(\mu)$ is equal to $\pm 2$.
We show that $\rho(\gamma)\ne \pm I$. There are two cases.\\
Case I. $\rho(\mu)=\pm I$. We know that the $1/0$ Dehn surgery
produces $S^3$. Hence $\rho$ induced a representation of
$\pi_1(S^3)$ in $PSL_2(\mathbb{C})$. It is trivial. Therefore, the
image of $\rho$ in $SL_2(\mathbb{C})$ is contained in $\{\pm I\}$.
So $\rho$ is reducible. This contradiction shows that Case I can
not happen.\\
Case II. $\rho(\mu)$ is parabolic. Since $\rho(\mu)$ and
$\rho(\lambda)$ commute, up to conjugation, we can assume that
\begin{equation*}
\rho(\mu)=\pm
  \left[
  \begin{matrix}
    1 & 1\\
    0 & 1
  \end{matrix}
  \right], \;\;
\rho(\lambda)=\pm
  \left[
  \begin{matrix}
    1 & t(\rho)\\
    0 & 1
  \end{matrix}
  \right].
\end{equation*}

Now $\gamma=p\mu+q\lambda$, so
$\rho(\gamma)=\rho(\mu)^p\rho(\lambda)^q$. If $t(\rho)=0$,
$\rho(\gamma)\ne \pm I$ unless $p=0$. If $t(\rho)\ne 0$, by the
assumption $t(\rho)\notin \mathbb{Q}$, then $\rho(\gamma)\ne \pm I$.
\end{proof}

\begin{proposition}\label{redu}
Suppose $\chi\in Z_{\gamma}$ is the character of a reducible
representation $\rho$. Then $\chi(\mu)\ne \pm 2$ and
$\rho(\gamma)=\pm I$.
\end{proposition}

\begin{proof}
Let $m$ be an eigenvalue of $\rho(\mu)$. By
\cite[Proposition~6.2]{CCGLS}, $m^2$ must be a root of the Alexander
polynomial $\Delta(t)$ of the knot. It is well-known that
$\Delta(1)\ne 0$, hence $m\ne \pm 1$ and $\chi(\mu)\ne \pm 2$. This
also means that $\rho(\mu)$ is not parabolic or $\pm I$. Since
$\rho(\mu)$ and $\rho(\gamma)$ commute and the trace of
$\rho(\gamma)$ is $\pm 2$, we must have $\rho(\gamma)=\pm I$.
\end{proof}

\begin{remark}
Compare with Proposition \ref{kin}, for the reducible characters, it
is much simpler and there is no need of the conjecture.
\end{remark}

Let $E(p,q)$ be the reducible curve $m^pl^q=\pm 1$ in
$\mathbb{C}^{*}\times \mathbb{C}^{*}$ for $p,q$ coprime integers.
Then the image $t_D(E(p,q))$ is a curve in $X(\partial M_K)$. We
know $r(X_0)$ is an irreducible curve in $X(\partial M_K)$. They do
not have common irreducible component because the traces of
characters of $X_0$ are not constant. Hence $t_D(E(p,q))\cap r(X_0)$
is finite. The set $t_D(E(p,q))\cap r(X_0)$ consists of possible
characters in $X_0$ which can also be the characters of $K(p/q)$,
where $K(p/q)$ denotes the closed $3$-manifold obtained from $M_K$
by the Dehn surgery along the simple closed curve of $\partial M_K$
which represents the class $\gamma$ in $H_1(\partial
M_K;\mathbb{Z})$. The following definition should be thought of as
the \emph{algebro-geometric} invariant for the $(p,q)$ Dehn surgery
of $M_K$.

\begin{definition}
\[
  b(p,q)=\sum_{\chi\in t_D(E(p,q))\cap r(X_0)}n_{\chi},
\]
where $n_{\chi}$ is the intersection multiplicity at $\chi$.
\end{definition}

\begin{theorem}\label{b-inv}
The integer $b(p,q)$ is a well-defined invariant of the $3$-manifold
$K(p/q)$ resulting from the Dehn filling on the hyperbolic knot
complement $M_K$. It depends on the knot $K$ and the surgery
coefficient $p/q$, and it is always positive.
\end{theorem}

\begin{proof}
As mentioned above, the set $t_D(E(p,q))\cap r(X_0)$ is finite and
hence $b(p,q)<\infty$. On the other hand, because
$\chi_{\rho_0}(\mu)=\pm 2$, $\chi_{\rho_0}$ is always contained in
this set, where $\rho_0\in R_0$ is the discrete faithful
representation of the hyperbolic metric. Therefore, $b(p,q)>0$.

The intersection between $t_D(E(p,q))$ and $r(X_0)$ is taking place
in the surface $X(\partial M_K)$ described in Proposition
\ref{torus}. If two hyperbolic knots $K_1$ and $K_2$ are
homeomorphic, then they have the isomorphic fundamental groups of
their complements in $S^3$, hence they have isomorphic $X_0$.
Therefore, they have the same $b(p,q)$. Thus $b(p,q)$ is an
invariant of $K(p/q)$ depending only on the hyperbolic knot $K$ and
the Dehn surgery coefficient $p/q$.
\end{proof}

Now set $S(p,q)=\{\chi\in t_D(E(p,q))\cap r(X_0)|\chi(\mu)\ne \pm
2\}\subset r(X_0)\subset X(\partial M)$.

\begin{proposition}\label{count}
Suppose the conjecture 1 is true. Then $S(p,q)$ is exactly the set
of characters in $X_0$ which are also the characters of $K(p/q)$.
\end{proposition}

\begin{proof}
This follows from the definition of $S(p,q)$ and Proposition
\ref{kin}.
\end{proof}

\begin{definition}\label{lambda}
$$\lambda (p,q)=\sum_{\chi\in S(p,q)} n_{\chi},$$
where $n_{\chi}$ is the intersection multiplicity at $\chi$.
\end{definition}

By Proposition \ref{redu}, the set $S(p,q)$ contains all possible
reducible characters. Hence the number $\lambda (p,q)$ counts both
irreducible and reducible characters of $K(p/q)$. We think it is
important to count both irreducible and reducible representations in
$SL_2(\mathbb{C})$ case. It is easy to count the abelian and
non-abelian reducible characters, so there is a computable way to
count the irreducible characters for $K(p/q)$.

\begin{theorem}\label{c-inv}
Assume the conjecture 1 is true. The quantity $\lambda (p,q)$ is a
well-defined \emph{algebro-geometric} $SL_2(\mathbb{C})$ Casson type
invariant of $K(p/q)$.
\end{theorem}

\begin{proof}
It follows from Theorem \ref{b-inv} and Proposition \ref{count}.
\end{proof}

Note that by definition, $\lambda (p,q)\leq b(p,q)$ for any coprime
$p$, $q$ and a hyperbolic knot in $S^3$.

The invariant $\lambda (p,q)$ is for the hyperbolic knot and its
$(p,q)$ Dehn surgery. An $SL_2(\mathbb{C})$ knot invariant obtained
from the character variety of $1$-dimensional components is given by
the first author in \cite{Li}. The construction in \cite{Li} was
purely topological by choosing generic smooth perturbations and
generic values of $\chi(\mu)$. The topological definition of the
Casson-type invariant is not easy to calculate. Our invariant
$\lambda (p,q)$, defined via the intersection multiplicity, is
easier or at least very explicit for computation. It can also be
interpreted as the intersection cycles of the appropriate cohomology
classes of $X(\partial M)$.

In section $4$, for some coefficients $(p,q)$, we compute the
invariants $b(p,q)$ and $\lambda(p,q)$ for the figure-eight knot.
See the table after Corollary \ref{8reim}.

\begin{proposition}\label{prop-liv}
If $\lambda (p,q)>0$, then $\pi_1(K(p/q))$ is non-cyclic.
\end{proposition}

\begin{proof}
If $\lambda (p,q)>0$, then there exists $\rho\in R_0$ such that its
character $\chi_{\rho}\in S(p,q)$. Hence $f_{\mu}(\chi_{\rho})\ne
0$, and $f_{\gamma}(\chi_{\rho})=0$. By
\cite[Proposition~1.5.2]{CGLS}, there exists a representation from
$\pi_1(K(p/q))$ to $PSL_2(\mathbb{C})$ with non-cyclic image. The
result follows.
\end{proof}

\begin{remark}
(i) Our invariant $\lambda (p,q)$ is defined as the algebraic
intersection multiplicity in $X(\partial M)$ from the $X_0$
component. This is different from the definition in \cite{Cu1} where
the number $\lambda(K_{p/q})$ is defined over all components
of $X(M)$ and the intersection is taken in different space.\\
(ii) For hyperbolic knot $K$, $K(p/q)$ may not be NSL manifold. For
the component $X_0$, the intersection in \cite{Cu1} is $r(X_0)\cap
(t_D(t_D^{-1}(\overline{r(X_0)})\cap \{m^pl^q=1\}))$. This is
different from our Definition \ref{lambda}. Moreover, the
intersection multiplicity of an intersection point in $X(\partial
M_K)$ is counted by its multiplicity in $X(K(p/q))$ via an
appropriate Heegaard splitting of $K(p/q)$ in \cite{Cu1}. \\
(iii) It would be interesting to prove that the two character
varieties of Heegaard handle-bodies are smooth and their
intersection is always proper. In \cite[Section~4]{FM}, Fulton and
MacPherson only sketched an argument that if the two smooth
subvarieties intersect properly, then the topological intersection
multiplicity agrees with the algebraic intersection multiplicity.
\end{remark}

In \cite[Section~1.4]{CGLS}, a norm $|.|$ is defined on the real
vector space $H_1(\partial M_K,\mathbb{R})$ with the property that
$|\gamma|=\text{deg}f_{\gamma}$ for any $\gamma\in H_1(\partial
M_K,\mathbb{Z})$. This norm is called the Culler-Shalen norm. In
particular, $|\gamma|=\text{deg}f_{\gamma}^{\prime}$.

Let $\pi: \widetilde{Y_0}\rightarrow Y_0$ be the birational
isomorphism. Note that $\pi$ is well-defined only on a Zariski dense
subset of $\widetilde{Y_0}$ and is surjective. A point of the set
$I=\widetilde{Y_0}\setminus  \pi^{-1}(Y_0)$ is called an ideal
point. Denote by $Z_{\gamma}^{\prime}$ the set of zeroes of the
meromorphic function $f_{\gamma}^{\prime}:\widetilde{Y_0}\rightarrow
\mathbb{C}$. Then
\[
  \text{deg}f_{\gamma}^{\prime}=\sum_{y\in
  Z_{\gamma}^{\prime}}v_y,
\]
where $v_y$ is the order of vanishing of $f_{\gamma}^{\prime}$ at
$y$.\\

Set $Z_1=\{y\in Z_{\gamma}^{\prime}|\pi(y)\in S(p,q)\}$ and
$I(p,q)=Z_{\gamma}^{\prime}\cap I$. Now we have another quantity:
\begin{equation}\label{v-num}
  \widehat{\lambda}(p,q)=\sum_{y\in Z_1}v_y.
\end{equation}
The natural question is to find the relationship between
$\widehat{\lambda}(p,q)$ and $\lambda(p,q)$ of the Definition
\ref{lambda}. It seems that there is no easy answer to this
question. See the remark below. Nevertheless, we have the following:
\begin{proposition}\label{deg-ine}
Assume that for every $\chi\in S(p,q)$, $t_D(E(p,q))$ intersects
$r(X_0)$ transversely at $\chi$. Then\\
(i) $\lambda(p,q)\leq \widehat{\lambda}(p,q)$,\\
(ii) $\lambda(p,q)+\widehat{I}(p,q)\leq \text{deg}f_{\gamma}$, where
$\widehat{I}(p,q)=\sum_{x\in I(p,q)}v_x$.
\end{proposition}

\begin{proof}
If the intersection is transverse, then intersection multiplicity
$n_{\chi}=1$. The order of vanishing is at least one,  hence (i)
holds.

For (ii), note that $Z_1$ and $I(p,q)$ are subsets of
$Z_{\gamma}^{\prime}$. Thus, $\lambda(p,q)+\widehat{I}(p,q)\leq
\text{deg}f_{\gamma}^{\prime}=\text{deg}f_{\gamma}$.
\end{proof}

\begin{remark}
If the intersection is not transversal, for instance, some $\chi$ is
a singular point of $r(X_0)$, then we do not know how to compare
them. One difficulty is that when $\chi$ is a singular point, there
is NO well-defined notion of the order of vanishing of
$f_{\gamma}^{\prime}$ at $\chi$ in $Y_0$. Moreover, $\pi^{-1}(\chi)$
has more than one element in $\widetilde{Y_0}$ and each of them is a
zero of $f_{\gamma}^{\prime}$.
\end{remark}

\section{An Example: The Figure Eight Knot} Throughout this
section, we shall denote by $M$ the complement of the figure-eight
knot in $S^3$. We compute its character variety and give some
applications.

It is well-known that $\pi_1(M)$ is given by two generators and one
relation:
\begin{equation}\label{8gp}
  \pi_1(M)=<\alpha,\beta|R(\alpha,\beta)>,\;
  R(\alpha,\beta)=\beta^{-1}\alpha^{-1}\beta\alpha\beta^{-1}\alpha\beta\alpha^{-1}\beta^{-1}\alpha,
\end{equation}

where $\alpha$, $\beta$ are meridians, and they are conjugate to
each other $\beta=\delta\alpha\delta^{-1}$ with
$\delta=\alpha^{-1}\beta\alpha\beta^{-1}$.

Set $\tau=\alpha\beta^{-1}\alpha^{-1}\beta$, we have a peripheral
subgroup $\pi_1(\partial M)$
\begin{equation}\label{8bgp}
 \pi_1(\partial M)=<\alpha,\lambda>; \;
 \lambda=\tau^{-1}\delta=\beta^{-1}\alpha\beta\alpha^{-1}\alpha^{-1}\beta\alpha\beta^{-1},
\end{equation}

where $\alpha$ is the meridian and $\lambda$ is the longitude. Note
in this section, we use $\alpha$ for the meridian instead of $\mu$.

Let us consider its character variety $X(M)$. We use (\ref{8gp}) for
the presentation of $\pi_1(M)$. For a
$SL_2(\mathbb{C})$-representation $\rho\in R(M)$, its character
$\chi$ is determined by $\chi(\alpha)$, $\chi(\beta)$, and
$\chi(\alpha\beta)$. We have the morphism $t:R(M)\rightarrow
\mathbb{C}^3$, $t(\rho)=(\chi(\alpha),
\chi(\beta),\chi(\alpha\beta))$ and $X(M)$ is the image $t(R(M))$.
Let $(x,y,z)$ be the affine coordinate for $\mathbb{C}^3$.

\begin{proposition}\label{8cha}
The affine variety $X(M)\subset \mathbb{C}^3$ is defined by the
following equations:
\begin{equation}\label{8d1}
  x=y,
\end{equation}
\begin{equation}\label{8d2}
 (x^2-z-2)(z^2-(1+x^2)z+2x^2-1)=0.
\end{equation}
In particular, we can identify $X(M)$ with the affine plane curve
$\{(x,z)\in \mathbb{C}^2|(x^2-z-2)(z^2-(1+x^2)z+2x^2-1)=0\}$. It has
two irreducible components.
\end{proposition}

\begin{proof}
Since $\alpha$ is conjugate to $\beta$,
$x=\chi(\alpha)=\chi(\beta)=y$. Hence (\ref{8d1}) follows. For the
second equation, by \cite[Theorem~1]{Wh}, the factor $x^2-z-2$
corresponds to characters of abelian representations, and the other
one corresponds to the characters of non-abelian representations.
\end{proof}

By the preceding proposition, the component defined by the equation
$z^2-(1+x^2)z+2x^2-1=0$ consists of the characters of non-abelian
representations, in particular, it contains the discrete faithful
representations of the complete hyperbolic metric of $M$. We denote
this component $X_0$. Hence we have :

\begin{corollary}\label{8-cor1}
The component $X_0$ is an irreducible smooth curve in $\mathbb{C}^2$
with the defining equation:
\begin{equation}\label{8d3}
   z^2-(1+x^2)z+2x^2-1=0.
\end{equation}
\end{corollary}

\begin{proof}
We only need to check that (\ref{8d3}) defines a smooth curve. Let
$f(x,z)=z^2-(1+x^2)z+2x^2-1$. Then $\frac{\partial f}{\partial
x}=-2xz+4x$ and $\frac{\partial f}{\partial z}=2z-1-x^2$. It is
straightforward to check that their is no common solution to the
equations $\frac{\partial f}{\partial x}=\frac{\partial f}{\partial
z}=f(x,z)=0$. Hence, the curve is smooth.
\end{proof}

\begin{corollary}\label{8-cor2}
There are exactly two reducible characters on $X_0$ and they
correspond to the points $(\pm \sqrt{5},3)$.
\end{corollary}

\begin{proof}
Let $\rho$ be a reducible representation whose character lies on
$X_0$. Let $m$ be the eigenvalue of $\rho(\alpha)$. By
\cite[Proposition~6.2]{CCGLS}, $m^2$ must be a root of the Alexander
polynomial $\Delta(t)=t^2-3t+1$ of the figure-eight knot. Thus,
\[
  m=\pm \sqrt{\frac{3\pm \sqrt{5}}{2}}=\frac{\pm \sqrt{5}\pm
  1}{2}; \, \text{and}\; x=\sigma(\rho(\alpha))=m+m^{-1}=\pm \sqrt{5}.
\]

Plug in (\ref{8d3}), we get $z=3$ with multiplicity two. Now it is
easy to check that $(\pm \sqrt{5},3)$ are exactly the intersection
points of $X_0$ and the component $x^2-z-2=0$. Since the latter
component consists of abelian characters, the result follows.
\end{proof}

\begin{remark}
Two inequivalent reducible representations may have the same
character. The points $(\pm \sqrt{5},3)$ above are such examples.
They are the characters of both some abelian representation and
non-abelian reducible representation which are clearly not
equivalent.
\end{remark}

\begin{proposition}\label{8cha-1}
Suppose that $\chi\in X_0$ and $\chi(\alpha)=\pm 2$. Then $\chi$ is
the character of a discrete faithful representation.
\end{proposition}

\begin{proof}
Plug $x=\pm 2$ in (\ref{8d3}), we get $z^2-5z+7=0$ and hence
$z=\frac{5 \pm \sqrt{-3}}{2}$. The following representation
$\theta:\pi_1(M)\rightarrow PSL_2(\mathbb{C})$ was given in
\cite{Ril2}:
\begin{equation*}
\theta(\alpha)=\pm
  \left[
  \begin{matrix}
    1 & 1\\
    0 & 1
  \end{matrix}
  \right], \;\;
\theta(\beta)=\pm
  \left[
  \begin{matrix}
    1 & 0\\
    -\omega & 1
  \end{matrix}
  \right];
\end{equation*}
where $\omega=\frac{-1+\sqrt{-3}}{2}$ is a primitive cube root of
unity. By \cite[Theorem~1]{Ril2}, $\theta$ is an isomorphism. Thus,
$\theta$ is a discrete faithful representation of $\pi_1(M)$. Let
$\overline{\theta}$ be the complex conjugation of $\theta$. It is
straightforward to check that the four points $(\pm 2, \frac{5 \pm
\sqrt{-3}}{2})$ correspond to the $SL_2(\mathbb{C})$ characters of
the lifts of $\theta$ and $\overline{\theta}$. The result follows.
\end{proof}

Denote by $\mathbb{CP}^2$ the complex projective plane. We use
$[X:Y:Z]$ to represent its homogenous coordinates. We will identify
$\mathbb{C}^2$ with the open subset $\{[x:1:z]|(x,z)\in
\mathbb{C}^2\}$. Note this is different from the standard notation.

Let $\widetilde{X_0}$ be a smooth projective model of $X_0$. We have
the following explicit description of $\widetilde{X_0}$ in
$\mathbb{CP}^2$.
\begin{proposition}\label{8-smo}
$\widetilde{X_0}$ is an elliptic curve and defined by the equation
in $\mathbb{CP}^2$:
\begin{equation}\label{8sm}
   YZ^2-Y^2Z-X^2Z+2X^2Y-Y^3=0.
\end{equation}
\end{proposition}

\begin{proof}
Equation (\ref{8d3}) is the defining equation of $X_0$. Substitute
$x=\frac{X}{Y}$, $z=\frac{Z}{Y}$, and we get (\ref{8sm}). That is,
it is the projective closure of $X_0$ in $\mathbb{CP}^2$. It
suffices to show that it is smooth. Let $F(X,Y,Z)$ be the left-hand
side of (\ref{8sm}). It is elementary to check that except
$(0,0,0)$, there is no common solution to the equations
$\frac{\partial F}{\partial X}=\frac{\partial F}{\partial
Y}=\frac{\partial F}{\partial Z}=0$. Hence, it is smooth. Since
(\ref{8sm}) is a cubic equation in $\mathbb{CP}^2$,
$\widetilde{X_0}$ has genus $1$.
\end{proof}

\begin{corollary}\label{8-ide}
$X_0$ has exactly two ideal points $[1:0:0]$ and $[0:0:1]$.
\end{corollary}

\begin{proof}
The ideal points correspond to $Y=0$. This gives $X^2Z=0$. Thus,
$X=0$ or $Z=0$.
\end{proof}

 Next, we explicitly construct the restriction map
$r:X_0\rightarrow X(\partial M)$ induced by the inclusion
$i:\pi_1(\partial M)\rightarrow \pi_1(M)$.

\begin{theorem}\label{8res}
The map $r:X_0\rightarrow X(\partial M)\subset \mathbb{C}^3$ is
given by the formulas:
\[
  r(x,z)=(x,F(x,z),G(x,z));
\]
where $F(x,z)=x^4-5x^2+2$ is the trace of the longitude $\lambda$,
and $G(x,z)=(4x-x^3)z+(x^5-4x^3-x)$ is the trace of $\alpha
\lambda$.
\end{theorem}

We will postpone the proof of Theorem \ref{8res} to the end of this
section. Instead, we first discuss its applications. The reason is
that the proof itself is elementary and lengthy, and probably is
known to the experts. However, we can not find it in the literature.
So we include here for completeness.

\begin{corollary}\label{8reim}
 Let $(u,v,w)$ be the affine coordinates of $\mathbb{C}^3$. Then the
 image $r(X_0)$ in $X(\partial M)$ is closed, and it is defined by
 the following equations:
 \begin{equation}\label{8reim1}
  v=u^4 - 5u^2 + 2,
\end{equation}
\begin{equation}\label{8reim2}
  w^2-(u^5-5u^3+2u)w+(u^8-10u^6+29u^4-19u^2)=0.
\end{equation}
\end{corollary}

\begin{proof}
 Let $V=X(\partial M)\cap V_0$, where $V_0$ is the algebraic set
 in $\mathbb{C}^3$ defined by the equations (\ref{8reim1}) and
 (\ref{8reim2}). Then it suffices to show that $r(X_0)=V$.

By Theorem \ref{8res}, $r(X_0)=\{(u,v,w)\in X(\partial
M)|u=x,v=F(x,z),w=G(x,z)\}$. We need to eliminate the variables $x$
and $z$. Equation (\ref{8reim1}) clearly follows from the expression
of $F(x,z)$. From $G(x,z)$, we solve $z$ in terms of $u$ and $w$:
\[
  z=\frac{w-(u^5-4u^3-u)}{4u-u^3}.
\]
By Corollary \ref{8-cor1}, the defining equation of $X_0$ is
$z^2-(1+x^2)z+2x^2-1=0$. Therefore we have the equation:
\[
  [\frac{w-(u^5-4u^3-u)}{4u-u^3}]^2-(1+u^2)[\frac{w-(u^5-4u^3-u)}{4u-u^3}]+2u^2-1=0.
\]
Now get rid of the denominator and simplify the expression, we
obtain equation (\ref{8reim2}). But the issue is when the
denominator $4u-u^3=0$. For $x=u=0$, By Theorem \ref{8res}, the
corresponding point in $r(X_0)$ is $(0,2,0)$ and we can verify this
point lies in $V$. Similarly for $x=u=\pm 2$, the corresponding
points $(\pm 2,-2,\mp 2)$ lie in $V$. Hence $r(X_0)\subset V$. For
the other direction, let $(u,v,w)\in V$, if $u\ne 0,\pm 2$, then
take $x=u$, $z=\frac{w-(u^5-4u^3-u)}{4u-u^3}$, $r(x,z)=(u,v,w)$. For
$u=0$, we have $v=2$,$w=0$; for $u=\pm 2$, we have $v=-2$ and $w=\mp
2$. These points are in $r(X_0)$. Hence $r(X_0)=V$.
\end{proof}

Now we can use the computer software \verb"SINGULAR" to compute the
invariants $b(p,q)$ and $\lambda(p,q)$ for the figure-eight knot.
For more information about this software, we refer to \cite{GP} and
the web site \verb"http://www.singular.uni-kl.de". Use
\verb"SINGULAR" we can find the defining ideal of $E(p,q)$. The
defining equations of $r(X_0)$ is given by the Corollary
\ref{8reim}. Hence we can find their intersection with multiplicity.
The following list is some examples we compute using
\verb"SINGULAR".

\begin{center}
  \begin{tabular}{|r|r|r|}\hline
     $(p,q)$ &   $b(p,q)$ & $\lambda(p,q)$\\ \hline
      (2,1) & 8  &6\\ \hline
      (3,1) & 8  &6\\ \hline
      (4,1) & 6  &4\\ \hline
      (5,1) & 10 &8\\ \hline
      (6,1) & 12 &10\\ \hline
      (-1,2)& 16 &14\\ \hline
      (7,2) & 16 &16\\ \hline
      (8,3) & 16 &16\\ \hline
      (-9,4)& 24 &22\\ \hline
   \end{tabular}
 \newline
\end{center}

From the above table we see that $b(p, q)$ and $\lambda (p, q)$ may
be same for some $(p, q)$ surgeries along the figure eight knot. In
general the two algebro-geometric invariants $b(p, q)$ and $\lambda
(p, q)$ are not equal.

Next we compute the Culler-Shalen norm. Denote by $x_0=[0:0:1]$ and
$x_1=[1:0:0]$. By Corollary \ref{8-ide}, they are the only ideal
points of $X_0$.

\begin{lemma}\label{l-norm}
(1) Both $x_0$ and $x_1$ are poles of the meromorphic functions
$f_{\alpha}$ and $f_{\lambda}$. Their orders are:
$v_{x_0}(f_{\alpha})=v_{x_1}(f_{\alpha})=2$, and
$v_{x_0}(f_{\lambda})=v_{x_1}(f_{\lambda})=8$. \\
(2) The Culler-Shalen norm $|(1,0)|=4$ and $|(0,1)|=16$.
\end{lemma}

\begin{proof}
(1) By Theorem \ref{8res}, for $\chi=(x,z)\in X_0$,
$I_{\alpha}(x,z)=x$ and $I_{\lambda}(x,z)=F(x,z)=x^4-5x^2+2$. Since
$f_{\alpha}=I_{\alpha}^2-4$ and $f_{\lambda}=I_{\lambda}^2-4$, it is
sufficient to show that $v_{x_0}(I_{\alpha})=v_{x_1}(I_{\alpha})=1$,
and $v_{x_0}(I_{\lambda})=v_{x_1}(I_{\lambda})=4$.

On $\widetilde{X_0}$, substitute $x=\frac{X}{Y}$, $z=\frac{Z}{Y}$,
we get for $[X:Y:Z]\in \widetilde{X_0}$,
\[
   I_{\alpha}([X:Y:Z])=\frac{X}{Y},\, \text{and}\,\,
   I_{\lambda}([X:Y:Z])=\frac{X^4-5X^2Y^2+2Y^4}{Y^4}.
\]
First, we consider $x_0=[0:0:1]$. Let $U_3=\{[X:Y:Z]|Z\ne 0\}$. Then
it is an open subset containing $x_0$. $U_3$ is identified with
$\mathbb{C}^2$ via $x=\frac{X}{Z}$ and $y=\frac{Y}{Z}$, where
$(x,y)$ are affine coordinates of $\mathbb{C}^2$. Dividing by $Z^3$
the both sides of the equation (\ref{8sm}) of $\widetilde{X_0}$ and
substituting $x=\frac{X}{Z}$ and $y=\frac{Y}{Z}$, we get the affine
equation in $U_3=\mathbb{C}^2$:
\[
  g(x,y)=y-y^2-x^2+2x^2y-y^3=0.
\]
Under this identification, $x_0$ is the origin $(0,0)$,
$I_{\alpha}(x,y)=\frac{x}{y}$ and
$I_{\lambda}(x,y)=(\frac{x}{y})^4-5(\frac{x}{y})^2+2$. Since
$\frac{\partial g}{\partial y}=1-2y+2x^2-3y^2$ and $\frac{\partial
g}{\partial y}(0,0)\ne 0$, the function $x$ is a local parameter of
the local ring of regular functions at $(0,0)$. Solve $g(x,y)=0$, we
get
\[
  I_{\alpha}(x,y)=\frac{x}{y}=x^{-1}u(x,y),\, \text{where}\,
  u(x,y)=\frac{y^2+y-1}{2y-1},\, u(0,0)\ne 0;
\]
and
\[
  I_{\lambda}(x,y)=x^{-4}w(x,y),\, \text{where}\,
  w(x,y)=u^4-5x^2u^2+2x^4
  ,\, w(0,0)\ne 0.
\]
Therefore, $x_0$ is a pole of $I_{\alpha}$ of order $1$, and it is a
pole of $I_{\lambda}$ of order $4$.

For $x_1=[1:0:0]$. Let $U_1=\{[X:Y:Z]|X\ne 0\}$. $U_1$ is identified
with $\mathbb{C}^2$ via $y=\frac{Y}{X}$ and $z=\frac{Z}{X}$, where
$(y,z)$ are affine coordinates of $\mathbb{C}^2$. Similarly, we
obtain the affine equation in $U_1=\mathbb{C}^2$:
\[
  h(x,y)=yz^2-y^2z-z+2y-y^3=0.
\]
Now $x_0$ is the origin $(0,0)$, $I_{\alpha}(x,y)=y^{-1}$ and
$I_{\lambda}(x,y)=y^{-4}-5y^{-1}+2$. Since $\frac{\partial
h}{\partial z}(0,0)\ne 0$, the function $y$ is a local parameter of
the local ring of regular functions at $(0,0)$. Thus, $x_0$ is a
pole of $I_{\alpha}$ of order $1$, and it is a pole of $I_{\lambda}$
of order $4$.

(2) By definition,
$|(1,0)|=\text{deg}f_{\alpha}=v_{x_0}(f_{\alpha})+v_{x_1}(f_{\alpha})$,
by (1),  $|(1,0)|=2+2=4$. Similarly, $|(0,1)|=8+8=16$.
\end{proof}

Use this lemma, we can compute the Culler-Shalen norm $|(p,q)|$ for
any $\gamma=p\alpha+q\lambda\in H_1(\partial M,\mathbb{Z})$.

\begin{proposition}\label{p-8norm}
For each $\gamma=p\alpha+q\lambda\in H_1(\partial M,\mathbb{Z})$,
the Culler-Shalen norm $|\gamma|=|(p,q)|=2(|p+4q|+|p-4q|)$.
\end{proposition}

\begin{proof}
By Lemma \cite[~1.4.1, 1.4.2]{CGLS}, for each ideal point $x$, there
is a homomorphism $\phi_x:H_1(\partial M,\mathbb{Z})\rightarrow
\mathbb{Z}$, such that, for each $\gamma=p\alpha+q\lambda$,
\[
  v_x(f_{\gamma})=|\phi_x(\gamma)|,\, \text{and}\, |\gamma|=\sum_{x:\,\text{ideal
  point}}|\phi_x(\gamma)|;
\]
where $v_x(f_{\gamma})$ denotes the order of pole of $f_{\gamma}$ at
$x$. For our case, we have two ideal points $x_0$ and $x_1$. By
Lemma \ref{l-norm}, we get $|\phi_{x_0}(\alpha)|=2$,
$|\phi_{x_0}(\lambda)|=8$, $|\phi_{x_1}(\alpha)|=2$ and
$|\phi_{x_1}(\lambda)|=8$. Since $\phi_{x_i}$ are homomorphisms, we
obtain that $|\phi_{x_{i}}(\gamma)|=|\pm 2p \pm 8q|$, $i=1,2$.
Therefore, we obtain that either
\[
  |\phi_{x_{0}}(\gamma)|+|\phi_{x_{1}}(\gamma)|=|2p+8q|+|2p+8q|,
\]
or
\[
  |\phi_{x_{0}}(\gamma)|+|\phi_{x_{1}}(\gamma)|=|2p+8q|+|2p-8q|.
\]
We claim that the first case can not happen. Suppose not. Take
$\gamma=-4\alpha+\lambda$, then
$|\gamma|=|\phi_{x_{0}}(\gamma)|+|\phi_{x_{1}}(\gamma)|=0$. Hence
$\text{deg}f_{\gamma}=|\gamma|=0$. This is impossible because
$f_{\gamma}$ is not constant. Thus the claim is proved and
$|\gamma|=|(p,q)|=2(|p+4q|+|p-4q|)$.
\end{proof}

Next let $M(0)$ be the closed $3$-manifold obtained by the Dehn
surgery of $M$ along the longitude $\lambda$. It is known that
$M(0)$ admits an essential torus \cite[Page~200]{Boy}. By Theorem
\ref{8res}, we show that $M(0)$ has non-abelian infinite fundamental
group.
\begin{proposition}\label{0sur}
The fundamental group $\pi_1(M(0))$ is a non-abelian, infinite
group.
\end{proposition}

\begin{proof}
Suppose that $\rho\in R(M)$ with the property that its trace
$\sigma(\rho(\lambda))=2$. Now by Proposition \ref{8res}, we have:
\[
  x^4-5x^2+2=2
\]
where $x=\sigma(\rho(\alpha))$. Solve this equation, we get $x=0$ or
$x=\pm \sqrt{5}$. By Proposition \ref{8cha}, there exists $\rho_0\in
R(M)$, such that $\sigma(\rho_0(\alpha))=0$ and
$\sigma(\rho_0(\lambda))=2$. In particular the eigenvalues of
$\rho_0(\alpha)$ is $\pm i$. By \cite[Proposition~6.2]{CCGLS},
$\rho_0$ is an irreducible representation. On the other hand, since
$\rho_0(\alpha)$ and $\rho_0(\lambda)$ commute and $\rho_0(\alpha)$
is not parabolic, $\rho_0(\lambda)$ must be the identity matrix.
Hence $\rho_0$ induce a representation of $\pi_1(M(0))$. The
irreducibility of $\rho_0$ implies that $\pi_1(M(0))$ is not
abelian.

Similarly, let $\rho_{\sqrt{5}}\in R(M)$, such that
$\sigma(\rho_{\sqrt{5}}(\alpha))=\sqrt{5}$ and
$\sigma(\rho_{\sqrt{5}}(\lambda))=2$. By the same reason,
$\rho_{\sqrt{5}}(\lambda)$ equals the identity matrix. Hence
$\rho_{\sqrt{5}}$ induces a representation of $\pi_1(M(0))$. We can
check that the image of $\rho_{\sqrt{5}}$ in $SL_2(\mathbb{C})$ is
torsion-free. Therefore, $\pi_1(M(0))$ is not finite.
\end{proof}

\begin{remark} We know that $M(0)$ is not a hyperbolic
manifold. $0/1$ is the one of the ten exceptional surgery slopes of
the figure-eight knot. It is interesting to know that we can prove
its fundamental group is non-cyclic, non-abelian and infinite just
from the elementary computations.
\end{remark}

Let $M(3)$ be the closed $3$-manifold obtained by the Dehn surgery
of $M$ along the simple closed curve $\gamma=3\alpha+\lambda$.
\begin{lemma}\label{3/1-irre}
  $M(3)$ has exactly three irreducible $SL_2(\mathbb{C})$ characters.
\end{lemma}
\begin{proof}
 Since $\pi_1(M(3))=\langle \pi_1(M)|\alpha^3\lambda=1\rangle$,
 we have an embedding of character varieties
 $X(M(3))\hookrightarrow X(M)$. On the other hand, $X_0$ is the
 only component of $X(M)$ containing characters of non-abelian
 representations. Thus, the irreducible characters of $M(3)$ are
 contained in the set $S=\{\chi|\chi\in X_0,
 \chi(\alpha^3\lambda)=2\}$. By (\ref{tf3}), we have
 \begin{equation}
  \chi(\alpha^2\lambda)=\chi(\alpha)\chi(\alpha\lambda)-\chi(\lambda),
 \end{equation}
and
  \begin{equation}\label{3/1-t}
    \begin{aligned}
      \chi(\alpha^3\lambda)=&\chi(\alpha)\chi(\alpha^2\lambda)-\chi(\alpha\lambda)\\
                           =&(\chi(\alpha)^2-1)\chi(\alpha\lambda)-\chi(\alpha)\chi(\lambda).
    \end{aligned}
  \end{equation}
By Theorem \ref{8res}, we obtain
 \begin{equation}\label{3/1-eq1}
  (x^2-1)[(4x-x^3)z+x^5-4x^3-x]-x(x^4-5x^2+2)=2.
 \end{equation}
It is clear that the set $S$ is exactly the common solutions to
equations (\ref{3/1-eq1}) and (\ref{8d3}). Note when $x=1$,
(\ref{3/1-eq1}) holds and is independent of the values of $z$. From
(\ref{8d3}), when $x=1$, $z=1$. So $(1,1)\in S$.

We solve $z$ in terms of $x$ from (\ref{3/1-eq1}), then plug in
(\ref{8d3}) and simplify the expression, we get
\begin{equation}\label{3/1-eq2}
  C(x)=x^4-4x^3+2x^2+4x+1=(x^2-2x-1)^2=0.
\end{equation}
It has two solutions $x=1\pm \sqrt{2}$. Hence the set $S$ has three
elements. They are not equal to $\pm 2$ or $\pm \sqrt{5}$. By
Corollary \ref{8-cor2}, $x\ne \pm \sqrt{5}$ implies that each one is
an irreducible character. By Proposition \ref{8cha-1} and \ref{kin},
$x\ne \pm 2$ means that each one is also a character of $M(3)$. The
result follows.
\end{proof}

Let $sl_2(\mathbb{C})$ be the Lie algebra of $SL_2(\mathbb{C})$.
Then we have the adjoint representation $Ad:
SL_2(\mathbb{C})\rightarrow Aut(sl_2(\mathbb{C}))$. For a
representation $\rho:\pi_1(M(3))\rightarrow SL_2(\mathbb{C})$, let
$H^1(M(3); sl_2(\mathbb{C})_\rho)$ be the first cohomology group
with coefficients in $sl_2(\mathbb{C})$ twisted by the composition
$Ad\circ \rho$.

\begin{proposition}\label{3/1-prop1}
Suppose that $\rho:\pi_1(M(3))\rightarrow SL_2(\mathbb{C})$ is
irreducible. \\ Then $H^1(M(3);sl_2(\mathbb{C})_\rho)=0$.
\end{proposition}

\begin{proof}
By \cite[Theorem~1.1]{BW}, $M(3)$ is not toroidal. Since
$H_1(M(3);\mathbb{Z})$ is finite, $M(3)$ is a small Seifeit fibered
space. By the preceding Lemma \ref{3/1-irre}, $M(3)$ has irreducible
representations, so $\pi_1(M(3))$ is not cyclic. By
\cite[Proposition~7]{BZ1}, $H^1(M(3); sl_2(\mathbb{C})_\rho)=0$.
\end{proof}

\begin{remark} For the other eight exceptional surgery slopes $\pm 1$, $\pm 2$,
$\pm 3$, $\pm 4$, we also have the explicit irreducible
representations with infinite images. Hence their fundamental groups
are all non-abelian and infinite. We omit the details here. Hence,
Proposition \ref{3/1-prop1} holds also for $M(\pm1)$, $M(\pm2)$ and
$M(-3)$.
\end{remark}

Now we turn to the proof of Theorem \ref{8res}. We need to find
explicit expressions of the traces $F(x,z)$, $G(x,z)$ of
$\rho(\lambda)$ and $\rho(\alpha\lambda)$ respectively in terms of
traces of $\rho(\alpha)$ and $\rho(\alpha\beta)$.

We need some preliminary lemmas. For $A$, $B$, $C\in
SL_2(\mathbb{C}))$, we have the following identities on their traces
\cite{Wh}:
\begin{equation}\label{tf1}
  \sigma(AB)=\sigma(BA);
\end{equation}
\begin{equation}\label{tf2}
 \sigma(A)=\sigma(A^{-1});
\end{equation}
\begin{equation}\label{tf3}
  \sigma(AB)=\sigma(A)\sigma(B)-\sigma(AB^{-1});
\end{equation}
\begin{equation}\label{tf4}
 \sigma(ABC)=\sigma(A)\sigma(BC)+\sigma(B)\sigma(AC)+\sigma(C)\sigma(AB)-\sigma(A)\sigma(B)\sigma(C)-\sigma(ACB);
\end{equation}
\begin{equation}\label{tf5}
 \text{For}\, m\geq 2,
 \sigma(A^m)=\sigma(A^{m-1})\sigma(A)-\sigma(A^{m-2}).
\end{equation}

Notice that (\ref{tf1}) is true for any $n\times n$ matrices. It is
easy to see that (\ref{tf4}) and (\ref{tf5}) follow from
(\ref{tf3}). Equations (\ref{tf2}) and (\ref{tf3}) can be checked by
direct
computations.\\

Now for $\rho\in R(M)$, set $a=\rho(\alpha)$, $b=\rho(\beta)$,
$x=\sigma(a)=\sigma(b)$ and
$z=\sigma(ab)=\sigma(\rho(\alpha\beta))$.

\begin{lemma}\label{tr1}
(i) $\sigma(a^2)=x^2-2$.\\
(ii)
$\sigma(a^{-1}b)=\sigma(ba^{-1})=\sigma(b^{-1}a)=\sigma(ab^{-1})=x^2-z$.\\
(iii)
$\sigma(b^{-1}aba^{-1})=\sigma(a^{-1}bab^{-1})=z^2-x^2z+2x^2-2$.
\end{lemma}

\begin{proof}
(i) By (\ref{tf5}), $\sigma(a^2)=\sigma(a)^2-\sigma(I)$, where $I$
is the $2\times 2$ identity matrix.

(ii) By (\ref{tf1}), $\sigma(a^{-1}b)=\sigma(ba^{-1})$ and
$\sigma(b^{-1}a)=\sigma(ab^{-1})$ ; by (\ref{tf2}),
$\sigma(ab^{-1})=\sigma(ba^{-1})$; by (\ref{tf3}),
$\sigma(ab^{-1})=x^2-z$.

(iii) By the subdivision $(b^{-1}a)ba^{-1}$ and (\ref{tf4}), we have
\[
  \sigma(b^{-1}aba^{-1})=\sigma(b^{-1}a)\sigma(ba^{-1})+\sigma(b)\sigma(b^{-1})
  +\sigma(a^{-1})\sigma(b^{-1}ab)-\sigma(b^{-1}a)\sigma(b)\sigma(a^{-1})-\sigma(I).
\]
Therefore,
\[
  \sigma(b^{-1}aba^{-1})=(x^2-z)^2+x^2+x^2-(x^2-z)x^2-2=z^2-x^2z+2x^2-2.
\]
The proof for $\sigma(a^{-1}bab^{-1})$ is the same by the
subdivision $(a^{-1}b)ab^{-1}$, and we omit it.
\end{proof}

\begin{lemma}\label{tr2}
(i) $\sigma(a^{-1}a^{-1}b)=x(x^2-z)-x$;\\
(ii) $\sigma(ba^{-1}a^{-1}b)=x^4-zx^2-2x^2+2$;\\
(iii) $\sigma(ba^{-1}a^{-1}bab^{-1})=x^2-z$.
\end{lemma}

\begin{proof}
(i) By (\ref{tf3}), we have
\[
  \sigma(a^{-1}(a^{-1}b))=\sigma(a^{-1})\sigma(a^{-1}b)-\sigma(a^{-1}b^{-1}a)
\]
By Lemma \ref{tr1}, each term of the right-hand side is known.
Hence,
\[
  \sigma(a^{-1}a^{-1}b)=x(x^2-z)-x.
\]

(ii) By (\ref{tf3}),
\[
  \sigma(b(a^{-1}a^{-1}b))=\sigma(b)\sigma(a^{-1}a^{-1}b)-\sigma(bb^{-1}aa)=\sigma(b)\sigma(a^{-1}a^{-1}b)-\sigma(a^2)
\]
The formula follows.

(iii) By (\ref{tf4}), we have
\begin{align*}
  \begin{split}
    \sigma((ba^{-1})(a^{-1}b)(ab^{-1}))&=\sigma(ba^{-1})\sigma(a^{-1}bab^{-1})+\sigma(a^{-1}b)\sigma(1)+\sigma(ab^{-1})\sigma(ba^{-1}a^{-1}b)-\\
                                       & \sigma(ba^{-1})\sigma(a^{-1}b)\sigma(ab^{-1})-\sigma(a^{-1}b)
  \end{split}
\end{align*}

Plug in what we know on the right-hand side and simplify, we obtain
the formula.
\end{proof}

\begin{lemma}\label{tr3}
(i) $\sigma(ab^{-1}a)=\sigma(aab^{-1})=x^3-zx-x$;\\
(ii) $\sigma(aab)=xz-x$;\\
(iii) $\sigma(ab^{-1}aba^{-1})=x$
\end{lemma}

\begin{proof}
(i) $\sigma(a(b^{-1}a))=\sigma(a)\sigma(b^{-1}a)-\sigma(b)$, and $\sigma(a(ab^{-1}))=\sigma(a)\sigma(ab^{-1})-\sigma(aba^{-1})$;\\
(ii) $\sigma(a(ab))=\sigma(a)\sigma(ab)-\sigma(ab^{-1}a^{-1})$\\
(iii)
$\sigma(a(b^{-1}aba^{-1}))=\sigma((b^{-1}aba^{-1})a)=\sigma(b^{-1}ab)=x$.
\end{proof}

\begin{proof}[Proof of Theorem~\ref{8res}]
First, let us compute $F(x,z)=\sigma(\rho(\lambda))$, the trace of
$\rho(\lambda)$. By (\ref{tf4}), we have
\begin{align*}
  \begin{split}
     \sigma((b^{-1}a)(ba^{-1})(a^{-1}bab^{-1}))&=\sigma(b^{-1}a)\sigma(ba^{-1}a^{-1}bab^{-1})+\sigma(ba^{-1})\sigma(ab^{-1})+\\
                                               & \sigma(a^{-1}bab^{-1})\sigma(b^{-1}aba^{-1})
                   -\sigma(b^{-1}a)\sigma(ba^{-1})\sigma(a^{-1}bab^{-1})-\sigma(I)
  \end{split}
\end{align*}

By Lemmas \ref{tr1} and \ref{tr2}, we know all the terms on the
right-hand side. Notice that $(x,z)\in X_0$, hence
$z^2-(1+x^2)z+2x^2-1=0$. Now divide the right-hand side by
$z^2-(1+x^2)z+2x^2-1$, the remainder is $F(x,z)$. we calculate that $F(x,z)=x^4-5x^2+2$.\\

For the trace of $\rho(\alpha\lambda)$, we have
\begin{align*}
  \begin{split}
     \sigma((ab^{-1}a)(ba^{-1})(a^{-1}bab^{-1}))&=\sigma(ab^{-1}a)\sigma(ba^{-1}a^{-1}bab^{-1})+\sigma(ba^{-1})\sigma(aab^{-1})+\\
                                                & \sigma(a^{-1}bab^{-1})\sigma(ab^{-1}aba^{-1})-\sigma(ab^{-1}a)\sigma(ba^{-1})\sigma(a^{-1}bab^{-1})
                                                -\sigma(a).
  \end{split}
\end{align*}

By Lemmas  \ref{tr1}, \ref{tr2} and \ref{tr3}, we have all the terms
on the right-hand side. Then we divide the result of the right-hand
side by the polynomial $z^2-(1+x^2)z+2x^2-1=0$ and $G(x,z)$ equals
the remainder. We calculate that it is $(4x-x^3)z+(x^5-4x^3-x)$.
\end{proof}

\end{document}